\documentclass[12pt]{article}
\usepackage{amssymb}
\usepackage{latexsym}
\usepackage[all]{xy}
\newcounter{mytab}
\newenvironment{mytab}[1]{
\refstepcounter{mytab}
\begin{center}
Table \arabic{mytab}. #1
\end{center}
}{}
\newtheorem{THM}{Theorem}[section]
\newtheorem{PROP}{Proposition}[section]
\newtheorem{LEM}{Lemma}[section]
\begin{document}
%
%%%%%%%%%%%%%%%%%%%%%%%%%%%%%%%%%%%%%%%%%%%%%%%%%%%%%%%%%%%%%%%%
%  

\title{Reflection Groups and Polytopes over Finite Fields, III}
 
\author{B. Monson\thanks{Supported by NSERC of Canada Grant \# 4818}\\
University of New Brunswick\\
Fredericton, New Brunswick, Canada E3B 5A3
\and and\\[.05in]
Egon Schulte\thanks{Supported by NSA-grants H98230-05-1-0027 
and H98230-07-1-0005}\\
Northeastern University\\
Boston, Massachussetts,  USA, 02115}

\date{ \today }
\maketitle

\begin{center}
{\em With best wishes for our friend and colleague, J\"org Wills}
\end{center}

\begin{abstract}
When the standard representation of a 
crystallographic Coxeter group $\Gamma$ is reduced modulo an odd prime $p$, one 
obtains a finite group $G^p$ acting on some orthogonal space over $\mathbb{Z}_p$. 
If $\Gamma$ has  a string diagram,  then $G^p$ will often be the automorphism group
of a finite abstract regular polytope. In parts I and II we established the basics 
of this construction and enumerated the polytopes associated to   groups of rank at
most $4$, as well as all  groups of spherical or Euclidean type.
Here we extend the range of our earlier criteria for the
polytopality of $G^p$. Building on this we investigate the class of 
`$3$--infinity' groups of general rank, and  then complete  
a survey of those  locally toroidal polytopes
which can be described by our construction.

\medskip



\end{abstract}

\section{Introduction and Notation}
\label{intro}

The regular polytopes continue to be a rich source of beautiful
mathematical ideas.  Their combinatorial features, for instance,
have been  generalized in the theory of 
\textit{abstract regular polytopes}. Here we conclude a series of three papers
concerning the properties of (abstract) regular polytopes, as constructed from
orthogonal groups over finite fields. 
Our main goal  
is to complete a description of
the \textit{locally toroidal polytopes} provided by our construction 
(see Section~\ref{loctor}).
To that end, in Section~\ref{interprop}  
 we establish some new structural theorems concerning the
`polytopality' of orthogonal groups.   As a test case, 
we also apply our methods in Section~\ref{higher}
to an interesting family of polytopes of general rank $n$.

Let us begin with a review of the basic set up and key results from parts I and II
(\cite{monsch1} and \cite{monsch2}, respectively).
In \cite{monsch1},  we first surveyed some
of the essential properties of an abstract regular polytope $\mathcal{P}$,
referring to \cite{arp} for details. Crucially, for each such $\mathcal{P}$ the 
automorphism group $\Gamma(\mathcal{P})$ is equipped with a natural list of 
involutory generators and is further a very special quotient of a certain Coxeter
group $G$. (We say that $\Gamma(\mathcal{P})$ is a \textit{string C-group}.)
Since $\mathcal{P}$ can be uniquely
reconstructed from  $\Gamma(\mathcal{P})$,  we may therefore shift our focus.

Throughout,  then, $G = \langle r_{0},\ldots,r_{n-1} \rangle$ will be a 
possibly infinite,
crystallographic Coxeter group $[p_{1},p_{2},\ldots,p_{n-1}]$ with a string Coxeter 
diagram $\Delta_{c}(G)$ (with branches labeled $p_{1},p_{2},\ldots,p_{n-1}$, 
respectively), obtained from the corresponding abstract Coxeter group 
$\Gamma = \langle \rho_{0},\ldots,\rho_{n-1} \rangle$ via the standard 
representation on real $n$-space~$V$. (Very often $G$ will be infinite.)
For any odd prime $p$, we may reduce $G$ modulo $p$ to obtain a subgroup 
$G^p$ of $GL_n(\mathbb{Z}_p)$ generated by the modular images of the $r_i$'s. 
% here
We shall abuse notation by referring to the modular images of objects by 
the same name (such as $r_i$, $b_i$, $B =[b_i \cdot b_j]$, $V$, etc.). In particular, 
$\{b_i\}$ will denote the standard basis for $V = \mathbb{Z}_p^n$.  In any event, 
$G^p$ is a subgroup of the orthogonal group $O( \mathbb{Z}_p^n)$ of isometries for 
the (possibly singular) symmetric bilinear form $x \cdot y$,  the latter 
being defined on $\mathbb{Z}_p^n$ by means of the Gram matrix $B$.  
Likewise, each $r_{i}$ remains a reflection, although
we may write
$$ r_i(x) = x - 2 \, \frac{x \cdot b_i}{b_i \cdot b_i}\, b_i$$
only if $b_i^2 := b_i \cdot b_i  \not\equiv 0 \bmod\,p$.
Concerning this situation, we now
make a convenient definition: if $p\geq 5$, or $p=3$ but no branch of 
$\Delta_{c}(G)$ is marked $6$, then we say that $p$ is \textit{generic} 
for $G$. Indeed, in such cases, no node label $b_i^2$ of the diagram $\Delta(G)$ 
(for a basic system) is zero $\bmod\, p$,
and the corresponding \textit{root} $b_i$ is anisotropic. Also,
a change in the underlying basic system for $G$ has the effect of merely conjugating $G^p$ in $GL_n(\mathbb{Z}_p)$. On the other hand, in the non-generic case, in which $p=3$ and $\Delta_{c}(G)$ has some branch marked $6$, the group $G^p$ may depend essentially on the actual diagram $\Delta(G)$ taken for the reduction mod $p$. (Note that $p$ generic does not necessarily mean that $p \nmid |G|$,  or that certain subspaces of $V$ are non-singular, etc.) 

Now we confront two questions: what exactly is the finite reflection group $G^p$  and when
is it a string $C$-group (i.e. the automorphism group of a finite, 
abstract regular $n$-polytope $\mathcal{P} = \mathcal{P}(G^p)$)?
To help answer the first question, 
we recall from \cite[Thm. 3.1]{monsch1} that an irreducible group $G^p$ of the above sort, 
generated by $n \geq 3$ reflections,
must necessarily be one of the following:
\begin{itemize}
\item an orthogonal group $O(n,p,\epsilon) = O(V)$ or $O_j(n,p,\epsilon) = O_{j}(V)$, excluding the cases $O_1(3,3,0)$, $O_2(3,5,0)$, $O_2(5,3,0)$ (supposing for these three that 
$\mbox{disc}(V) \sim 1$), and also excluding the case $O_j(4,3,-1)$; or
\item the reduction mod $p$ of one of the finite linear Coxeter groups of type $A_n$ 
($p\nmid n+1$), $B_n$, $D_n$, $E_6$ ($p \neq 3$), $E_7$, $E_8$, $F_4$, $H_3$ or $H_4$.  
\end{itemize}
We shall say in  these two cases that $G^p$ is of \textit{orthogonal} or \textit{spherical type}, respectively, although there is some overlap for small primes.   
Concerning our groups $G^p$, it is only a slight abuse of notation to let 
$[p_1, \ldots , p_{n-1}]^p$ denote the modular representation of a 
group $[p_1, \ldots , p_{n-1}]$, so long as $p$ is generic for the group.

Let us turn to our second question.
The generators $r_i$ of $G^p$ certainly satisfy the Coxeter-type relations inherited from $G$. Thus 
$G^p$ is a \textit{string C-group}  if and only if it
satisfies  the following intersection property 
on standard subgroups:

\begin{equation}
\label{interII}
\langle r_i\,|\, i \in I \rangle \cap \langle r_i\,|\,i \in J \rangle =
\langle r_i\,|\, i \in I \cap J \rangle \;,
\end{equation}
for all $I, J \subseteq \{0, \ldots, n-1\}$ (see \cite[\S 2E]{arp}).
Our main problem is therefore to determine when $G^p$ satisfies
(\ref{interII}).  Before reviewing a few  preliminary results in this direction,
we establish some notation.

For any $J \subseteq \{0, \ldots, n-1\}$, we let 
$G_J^p := \langle r_j \,|\, j \not\in J \rangle $;  in particular, for 
$k,l \in \{0, \ldots n-1\}$ we let 
$G_k^p := \langle r_j \,|\, j \neq k\rangle$ and 
$G_{k,l}^p := \langle r_j \,|\, j \neq k,l \rangle$. 
 We also let $V_J$ be the subspace of 
$V = \mathbb{Z}_p^n$ spanned by $\{b_j \!\mid\! j \not\in J \}$, and 
similarly for $V_k, V_{k,l}$. Note that $V_J$ is $G_J^p$-invariant. In particular, 
$G_j^p$ acts on $V_j$, for $j=0$ or $n-1$. The upshot 
of Lemma 3.1 in \cite{monsch2} is that this action is faithful when $p$ is generic for $G$. 
Referring to \cite[Eq. 12]{monsch1}, we record here
a useful rule for inductively computing the determinant 
of the Gram matrix $B = [b_i \cdot b_j]$. Letting $B_J$ be the submatrix obtained
by deleting all rows and columns indexed by $J$, we have, for example,
\begin{equation}\label{discdet}
\det(B) = b_0^2 \det(B_0) - (b_0 \cdot b_1)^2 \det(B_{0,1})\; .
\end{equation}

We will frequently refer to  the following general properties of string $C$-groups,
here  as they apply to the groups   $G^p$:

\begin{PROP}\label{keyprops} 

{\rm (a)} $G^p$ is a string $C$-group if and only if
    $G^p_{0} , G^p_{n-1}$ are string $C$-groups and
    $G_0^p \cap G_{n-1}^p = G_{0,n-1}^p$.
    
{\rm (b)} If $G^p$ is a string $C$-group, then so too is any subgroup
    $G^p_J$, for $ J \subseteq \{0, \ldots , n-1\}$.
\end{PROP}
\noindent
\textbf{Proof}.  See \cite[2E16 and 2E12]{arp}. \hfill$\square$

In the next section we extract from \cite{monsch1} and \cite{monsch2} 
various  more specialized criteria for 
$G^p$ to be a string C-group. These concern the features of $V$
as an orthogonal space, as well as the action of standard
subgroups of $G^p$ on $V$. Using them, we were able in \cite{monsch1} to  classify 
all groups  $G^p$, and their polytopes, whenever $n \leq 3$, as
well as  when $G$ is of spherical or Euclidean type, for
all ranks $n$. Then in \cite{monsch2} we extended the classification to all cases in  rank $4$.
After generalizing these criteria, 
it will be clear that we have enough machinery to systematically 
extend our efforts  to polytopes of 
still higher rank. However, already in rank $4$ there is a bewildering variety of possibilities, 
so that
below we shall investigate only a few families of special interest.

\section{More on the Intersection Property}
\label{interprop}

%%%%%%%%%%%%%%%%%%%%

Let us review various situations  in which  $G^p$ is guaranteed to be a 
string $C$-group. First of all, this
will be the case if one of the subgroups $G_{0}^p$ or $G_{n-1}^p$ is 
spherical and the other is a string C-group:

\begin{THM}\label{intersphere} \mbox{\rm \cite[Th. 4.2]{monsch1}} 
Let $G = \langle r_0, \ldots , r_{n-1} \rangle$  
be a crystallographic linear Coxeter group with string diagram, 
and suppose the prime $p \geq 3$.
If $G_{n-1}$ is of spherical type
and $G_{0}^p$ is a string C-group,  or (dually) if
$G_0$ is of spherical type and $G_{n-1}^p$ is a string C-group,
then $G^p$ is a string C-group.
\end{THM}
We note that the proof supplied in \cite{monsch1} is  
inadequate for  the groups $G = [6,k]$ with $p=3$, though only for 
some of the possible basic systems (which need not be equivalent
in  these  non-generic cases). A familiar  example, taking $k=3$, is
the Euclidean group with diagram 
$$ \stackrel{1}{\bullet}\!\frac{}{\;\;\;\;\;\;}\!
	\stackrel{3}{\bullet}\!\frac{}{\;\;\;\;\;\;}\!
	\stackrel{3}{\bullet}\;. 
$$
Nevertheless, the intersection condition can be verified 
for all these groups, 
using  GAP \cite{gap}
or by hand.  In fact, such
peculiar exceptions  appear only peripherally in this 
paper.

The next two theorems utilize the occurence of groups of orthogonal type.

\begin{THM}\label{internonsing} \mbox{\rm \cite[Th. 4.1]{monsch1}} 
Let $G = \langle r_0, \ldots , r_{n-1} \rangle$  
be a crystallographic linear Coxeter group with string diagram, 
and suppose the prime $p \geq 3$. 
Suppose that $G_0^p$ and $ G_{n-1}^p$ are string C-groups, and that 
the subspace $V_{0,n-1}$ is non-singular.
Then if $G_{0,n-1}^p$ is the full orthogonal group $O(n-2,p,\epsilon)$
on $V_{0,n-1}$,  $G^p$ must be a string C-group. 

\end{THM}

%%%%%%%%%%%%%%%%%%%%
We now take a closer look at ways in which the geometry
of the various subspaces $V, V_0, V_{n-1}$ or $V_{0, n-1}$ affects the interaction
of the corresponding subgroups of $G^p$. 
The \textit{fully non-singular} case, proved in  \cite[Th. 3.2]{monsch2}
and generalized in part (a) below, 
sometimes allows
us to \textit{reject} large classes of 
groups $G^p$ as C-groups because of the size of their subgroups 
$G_{0}^p \cap G_{n-1}^p$.
When we leave the fully non-singular case, we  must adjust our approach in
various ways, depending on which of the various subspaces is
singular. The case in which just the \textit{middle section} is singular,
proved in  \cite[Th. 3.3]{monsch2} and repeated in (b) below, can sometimes
be used to affirm the polytopality of $G^p$ (see \cite[Cor. 3.2]{monsch2}).

In any ambient space $V$,  each non-singular subspace $W$ 
induces an orthogonal direct sum $V = W \perp W^\perp$ 
(\cite[Ch. 6, Lemma 2.1]{cohn}). Now consider
$O(W)$, the orthogonal group for $W$ (equipped
with the bilinear form inherited from $V$).
It is   easy to check that 
the mapping 
\begin{equation}
\label{kap}
\begin{array}{ccl}
\lambda: O(W)& \longrightarrow &  \mbox{Stab}_{O(V)} W^\perp \\
g  & \longrightarrow  & g \perp 1_{W^\perp} 
\end{array}
\end{equation}
establishes an isomorphism between $O(W)$ and a subgroup of
the \textit{pointwise}
stabilizer of $W^\perp$ in $O(V)$. 
We may therefore identify $O(W)$ with 
this subgroup; this is done without
much comment  for several subspaces $W$ in Theorem~\ref{subspcrit} below.
If $V$ happens to be non-singular, then
the spinor norm on $O(W)$ is also invariant under this 
identification \cite[Th. 5.13]{art}, and
we clearly have $O(W) \simeq \mbox{Stab}_{O(V)} W^\perp$. 

Let us now turn to the subgroup   
$O_1(W) := \langle r_a \, | \,  a \in W , a^2 = 1\rangle$.
By \cite[Prop. 3.1]{monsch1}, $O_1(W)$ almost always coincides 
with the kernel  of the spinor norm on $O(W)$ and so 
then has index $2$ in $O(W)$. However, for $\dim(W) \geq 2$ there
are two exceptions to this:
if $O(W)$ is isomorphic to either  
\begin{equation}\label{badgps}
[B_3]^3 \simeq O(3,3,0)\,\mbox{ (with disc $\sim 1$) or  } \;
[F_4]^3 \simeq O(4,3,+1)\; ,
\end{equation}
then $O_1(W)$ has index $3$ in the spinor kernel \cite[p. 301]{monsch1}.

In similar fashion we can work with a \textit{singular}
subspace $W$ of a non-singular space $V$. Here 
we let $\widehat{O}(W)$ denote the subgroup of $O(W)$ 
consisting of those isometries which act trivially on $\mbox{\rm rad}\, W$
(see \cite[Section 3]{monsch2}). 
It is not hard to show that $\widehat{O}(W)$  contains  and is generated by 
all reflections 
with non-isotropic roots in $W$. 
Furthermore,  we may  define a 
\textit{spinor norm} $\theta$ on $O(W)$; and 
$\widehat{O}_{1}(W)$  will   denote the subgroup 
of $\widehat{O}(W)$ generated by  
reflections in $O(W)$ with square  spinor norm.
In the proof of \cite[Th. 3.3]{monsch2}, we employed a 
variant of the mapping in (\ref{kap}) to show that
$\widehat{O}(W)$ can also be identified with a suitable subgroup  
of the pointwise stabilizer 
of $W^\perp$ in $O(V)$, so long as $W$ is a subspace $V_{0,n-1}$ 
(of codimension  $2$ in $V$); this is the only case that we require.
Again we find that $\widehat{O}_{1}(W)$ usually has index
$2$ in $\widehat{O}(W)$; for $\dim(W) \geq 2$, exceptions occur
when $O(W/{\mbox{rad}(W)})$ is either $O(2,3,+1)$ or one of the groups
in (\ref{badgps}).

\vspace*{5mm}
Let us  assemble our old results, along with some new criteria,
into one package:

\begin{THM}
\label{subspcrit}  
Let $G = \langle r_0, \ldots , r_{n-1} \rangle$  be a crystallographic linear 
Coxeter group with 
string diagram. Suppose  that $n\geq 3$,    that the prime $p$ is generic for $G$
and that there is a square among the labels of the nodes $1,\ldots,n-2$ of the 
diagram $\Delta (G)$ (this can be achieved by readjusting the node labels). 
For various subspaces $W$ of  $V$ we identify 
$O(W)$, $\widehat{O}(W)$, etc. with suitable subgroups of the  
pointwise stabilizer
of $W^\perp$ in $O(V)$.

\noindent{\rm (a)}
Let the subspaces  $V_{0}$, $V_{n-1}$ and  
$V_{0,n-1}$ be non-singular, and let $G_0^p$, $G_{n-1}^p$ be of orthogonal type.  
 
{\rm (i)}   Then $G_{0}^p \cap G_{n-1}^p$ acts trivially on $V_{0,n-1}^{\perp}$ and
$ O_{1}(V_{0,n-1}) \leq G_{0}^p \cap G_{n-1}^p \leq  O (V_{0,n-1})$. 

{\rm (ii)} If $G_{0}^{p} = O(V_{0})$ and $G_{n-1}^{p} = O(V_{n-1})$,  then 
$ G_{0}^p \cap G_{n-1}^p = O (V_{0,n-1})$. 

{\rm (iii)} If either $G_{0}^{p} = O_1(V_{0})$  or $G_{n-1}^{p} = O_1(V_{n-1})$,  then 
$ G_{0}^p \cap G_{n-1}^p = O_{1}(V_{0,n-1})$.   

\noindent{\rm (b)}
Let $V$, $V_{0}$, $V_{n-1}$ be non-singular, let $V_{0,n-1}$ be singular 
(so that $n \geq 4$), 
and let $G_0^p, G_{n-1}^p$ be of orthogonal type.  

{\rm (i)}  Then  $G_{0}^p \cap G_{n-1}^p$ acts trivially on $V_{0,n-1}^{\perp}$, and
$\widehat{O}_{1}(V_{0,n-1}) \leq G_{0}^p \cap G_{n-1}^p \leq  \widehat{O} (V_{0,n-1})$.

{\rm (ii)} If $G_{0}^{p} = O(V_{0})$ and $G_{n-1}^{p} = O(V_{n-1})$, then  
$ \widehat{O}(V_{0,n-1}) = G_{0}^p \cap G_{n-1}^p $. 

{\rm (iii)}  If either $G_{0}^{p} = O_1(V_{0})$ or  
$G_{n-1}^{p} = O_1(V_{n-1})$, then 
$ \widehat{O}_1(V_{0,n-1}) = G_{0}^p \cap G_{n-1}^p $.

\noindent{\rm (c)}  Suppose $V , V_{0,n-1}$ are non-singular 
while at least one of  $V_{0}$, $V_{n-1}$ is  singular.
Also suppose that $G_{0,n-1}^p$ is   of orthogonal type, with
 $G^p = O_1(V)$  when  $G_{0,n-1}^p = O_1(V_{0,n-1})$. Then 
$G_{0}^p \cap G_{n-1}^p  = G_{0,n-1}^p$.

\end{THM}
%%%%%%%%%%%%%%%%%%%%

\noindent\textbf{Proof}. When $V$ is non-singular,  parts (a)--(i),(ii) appear
as Theorem 3.2
in \cite{monsch2}. For $V$ singular, $\mbox{rad}(V) = \langle c \rangle$
is 1-dimensional, and we may choose a basis $w, w'$ for $V_{0,n-1}^\perp$
so that $c = w + w'$, $V_{n-1} = V_{0,n-1} \perp \langle w \rangle$,
$V_{0} = V_{0,n-1} \perp \langle w' \rangle$ and 
$V =  V_{n-1} \perp \langle v \rangle  = V_{0} \perp \langle v' \rangle$,
with $v = v' = c$. Then $g \in G_0^p \cap G_{n-1}^p$ implies
that $g(w) = \alpha w$, $g(w') = \alpha' w'$, where  $\alpha, \alpha' \in \{ \pm 1\}$.
Since $g(c) = c$, we have $ \alpha = \alpha' = 1$, so that
$g \in O(V_{0,n-1})$.  The rest of the proof of (i) and (ii)  
proceeds as in \cite[Th. 3.2]{monsch2}. 

For (a)--(iii) we may suppose
$G_{0}^{p} = O_1(V_{0})$. When $n=3$ there is nothing to prove, since
node $1$ has a square label and so $O(V_{0,2}) = O_1(V_{0,2})$. 
Now suppose $n-2 \geq 2$, so that there exists a reflection $r \in O(V_{0,n-1})$
with \textit{non-square} spinor norm. Since $r \not\in O_1(V_{0}) = G_{0}^{p}$,
we must by (i) have $ G_{0}^p \cap G_{n-1}^p = O_{1}(V_{0,n-1})$,
so long as $O_1(V_{0,n-1})$ has index $2$ in $O(V_{0,n-1})$.
As we observed earlier, this almost always holds. In fact, neither of the 
groups indicated  in (\ref{badgps}) can occur in our setup 
(as $O(V_{0,4}) , O(V_{0,5})$,
respectively).
Indeed, since $p=3$ in either case and since $G_{0}^{p} = G_0^3 = O_1(V_{0})$,
nodes $1, 2, \ldots , n-1$ must all be labelled by squares mod $3$. 
Also, $p$ is generic for $G$. These two restrictions imply that each of the standard
rotations $r_{j-1} r_j$ in $G$, except
possibly for $j=1$, must have period $3$ or $\infty$ 
(or $2$, if $\Delta(G)$ is disconnected).
In any case, $G_0^3 \simeq S_{a_1} \times \ldots \times  S_{a_k}$ is a direct
product
of $k \geq 1$ symmetric groups, where 
$ (a_1 -1) + \ldots + (a_k -1) = n-1$ ($ = 4$ or $5$ in the two cases).
A direct check of the possible orders shows that $G_0^3$ could not then 
be of orthogonal type in dimension $4$ or $5$ respectively.

In (b) we have $n \geq 4$; indeed, for $n=3$ we note that $V_{0,2}$ must be
non-singular, since $p$ is generic for $G$.  
Parts  (i),(ii) appear as Theorem 3.3 in \cite{monsch2}. 
We settle part  (iii)  
in much the same way as for (a)--(iii) above. Suppose that
$G_0^p = O_1(V_0)$ and let
$X := V_{0,n-1}/{\mbox{rad}(V_{0,n-1})}$, a non-singular space of 
dimension $n-3$. 
(Note that $X \simeq V_{0,1,n-1} \simeq V_{0,n-2,n-1}$.)
If $n=4$, the invariant quadratic form induced 
on $X$ must
be equivalent to $x_1^2$; then $\widehat{O}(V_{0,3}) = \widehat{O}_1(V_{0,3})$
and (b)--(iii) follows trivially. Now suppose that $n \geq 5$.
By our earlier remarks,  $\widehat{O}_1(V_{0,n-1})$ usually has index $2$ in 
$\widehat{O}(V_{0,n-1})$, in which case (b)--(iii) follows easily.
The three exceptional cases have $p=3$ with $n = 5,6,7$.
But as in part (a)--(iii) above, these
groups cannot occur as $O(X)$ when $G_0^p = O_1(V_0)$ or
 $G_{n-1}^p = O_1(V_{n-1})$ .

In part (c) there are two very similar cases, depending on whether one or both 
of $V_{0}$, $V_{n-1}$ are  singular.  To begin with, each
$g \in G_{0}^p \cap G_{n-1}^p$ certainly fixes ${\rm rad}(V_0)$ and
${\rm rad}(V_{n-1})$ pointwise. It is then easy to show in the two cases 
that $g$ fixes $V_{0 , n-1}^\perp$ pointwise, so that $g \in O(V_{0,n-1})$. 
Next one shows that
$ O_1(V) \cap O(V_{0,n-1}) = O_1(V_{0,n-1})$,
using  the assumption on square labels and \cite[Th. 5.13]{art}. 
(Here, too,  we must consider, and again exclude, the possibility 
that $O(V_{0,n-1})$  is one of the  
groups in (\ref{badgps}).)
Since either 
$G^p = O_1(V)$ or  $G_{0,n-1}^p = O(V_{0,n-1})$, we now have 
$g \in G_{0,n-1}^p$.
\hfill$\square$

%%%%%%%%%%%%%%%%%%%%
\vspace*{5mm} 
%%%%%%%%%%%%%%%%%%%%
Theorem~\ref{subspcrit} has several immediate and useful consequences. For example, in \cite[Cor. 3.2]{monsch2},  we used a preliminary version of part (b)(iii) to prove that
$[k,\infty,m]^p$ is a $C$-group for any odd prime $p$ and integers
$k,m \geq 2$. On the other hand, part (a)(iii) led  just as easily to a proof that $[\infty, 3, \infty]^p$ is a $C$-group only when $p = 3,5,7$.

Next we generalize \cite[Th. 3.4]{monsch2}, which concerns $4$-polytopes for which 
the facet group $G_3$ (say) is  Euclidean
and so situated   that the `point group' acts on the middle section of the polytope. Our first step is a closer look at the geometric action of groups
of affine Euclidean isometries.
In the background we typically have an abstract  Coxeter group of Euclidean (or `affine') type, faithfully represented in the standard way as a linear reflection group
$E = \langle r_0, \ldots , r_m \rangle$ on  real $(m+1)$-space
$W$.  Recall that $E$ preserves a positive semidefinite form $x \cdot y$, so that
$\mbox{rad}(W) = \langle c \rangle$ is $1$-dimensional. Since $r_j(c) = c$, for $0\leq j \leq m$, $E$ is in fact a subgroup of $\widehat{O}(W)$.

To actually exploit the structure of $E$ as a group of isometries on Euclidean $m$-space, we pass to the contragredient representation of $E$ in the dual space $\check{W}$ (as in \cite[5.13]{humph}). 
Since $c$ is fixed by $E$, we see that $E$ leaves invariant any translate of the $m$-space
$$ U =  \{ \mu \in \check{W} \,:\, \mu(c) =0 \}\;.$$
Next,  for each $w \in W$ define $\mu_w \in \check{W}$ by $\mu_w(x) := w \cdot x$. The mapping
$w \mapsto \mu_w$ factors to  a linear isomorphism between $W/\mbox{rad}(W)$ and $U$, and  so we transfer to $U$ the \textit{positive definite} form induced by  $W$ on $W/\mbox{rad}(W)$. 
Now choose any $\alpha \in \check{W}$ such that $\alpha(c) = 1$, and let
$\mathbb{A}^m := U +\alpha$. Putting all this together we may now
think of $\mathbb{A}^m$ as \textit{Euclidean $m$-space},
with $U$ as its \textit{space of translations}. Indeed, each fixed $\tau \in U$
defines an isometric translation on $\mathbb{A}^m$:
$$\mu \mapsto \mu +\tau,\;\; \forall \mu \in \mathbb{A}^m\;.$$
It is easy to check that this mapping on $\mathbb{A}^m$ is induced 
by a unique 
isometry $t \in \widehat{O}(W)$, namely the \textit{transvection}
\begin{eqnarray*}
t(x) & = &x - \tau(x) c, \\
     & = &x - (x\cdot a)\, c ,
\end{eqnarray*}
where $\tau = \mu_a$ for suitable $a \in W$. (Remember  here that we employ 
the contragredient representation of $\widehat{O}(W)$ on $\check{W}$, not just that of $E$.)
In summary, we can therefore safely
think of translations as transvections.

In the following table we list those Euclidean Coxeter groups which are  
relevant to our analysis (see \cite[\S 6B]{monsch1}). Concerning the 
group 
$E = [4, 3^{m-2}, 4]$ (for the familar cubical tessellation of 
$\mathbb{A}^m$), we recall our convention that
$3^{m-2}$ indicates a string of $m -2 \geq 0$ consecutive $3$'s.
\begin{center}
\begin{tabular}{c|c|c|c}\label{euctypes}
 The group $E$&  $m = \dim(\mathbb{A}^m)$& One possible diagram& 
 The corresponding vector\\
& & $\Delta(E)$ &  $c \in \mbox{rad}(W)$  \\\hline
 & &  & \\
 $[4, 3^{m-2}, 4]$ & $m \geq 2$&
$ \stackrel{2}{\bullet}\!\frac{}{\;\;\;\;\;\;}\!
	\stackrel{1}{\bullet}\!\frac{}{\;\;\;\;\;\;}\!
	\stackrel{1}{\bullet} \; \cdots \;
	\stackrel{1}{\bullet}\!\frac{}{\;\;\;\;\;\;}\!
	\stackrel{1}{\bullet}\!\frac{}{\;\;\;\;\;\;}\!
	\stackrel{2}{\bullet}$
                      & $c = b_0 + 2(b_1 + \ldots + b_{m-1}) + b_{m}$  \\ \hline
& & &\\
$[3,3,4,3]$        & $4$& 
	$ \stackrel{1}{\bullet}\!\frac{}{\;\;\;\;\;\;}\!
	\stackrel{1}{\bullet}\!\frac{}{\;\;\;\;\;\;}\!
	\stackrel{1}{\bullet}\!\frac{}{\;\;\;\;\;\;}\!
	\stackrel{2}{\bullet}\!\frac{}{\;\;\;\;\;\;}\!
	\stackrel{2}{\bullet}$
		&  $c = b_0 + 2 b_1 + 3 b_2 + 2 b_3 + b_4$  \\ \hline
& &  &\\
$[3,6]$            &$2$ &
$ \stackrel{1}{\bullet}\!\frac{}{\;\;\;\;\;\;}\!
	\stackrel{1}{\bullet}\!\frac{}{\;\;\;\;\;\;}\!
	\stackrel{3}{\bullet}$			 & $c = b_0 + 2 b_1 + b_2$ \\ \hline
& &  &\\
$[\infty]$ & $1$ & 
$\stackrel{1}{\bullet} =\!\!=\!\!= \stackrel{1}{\bullet}$ & $c = b_0 +  b_1$ \\ %\hline	 
\end{tabular}

\vspace*{4mm}
 
\end{center}

\begin{mytab}{Euclidean Coxeter Groups}\label{EucCoxGps}\end{mytab}

\noindent
An investigation of
the action of these discrete reflection groups on the Euclidean 
$m$-space $\mathbb{A}^m$ shows, in each case, that $E$ splits as the semidirect product of 
the (normal) subgroup $T$ of translations with a certain (finite) \textit{point group} group $H$:
\begin{equation}\label{eucsemi} 
E \simeq T \rtimes H\;. 
\end{equation}
(See \cite [Prop. 4.2]{humph}.) We can and do display each  group in the table so that  
 $H = E_0 = \langle r_1, \ldots, r_m \rangle$.

Returning now to our generalization of \cite[Th. 3.4]{monsch2}, we 
 suppose that $G_{n-1}$ is of Euclidean type.  
Of course, a dual result holds when $G_0$ is Euclidean.

\begin{THM}
\label{eucgen} 
Let $G = \langle r_0, \ldots , r_{n-1} \rangle$  be a crystallographic linear Coxeter group 
with string diagram. Suppose that $G_{n-1}$  is Euclidean, 
with $G_{n-1} = T \rtimes G_{0,n-1}$, where $T$ is the translation subgroup
of $G_{n-1}$. Suppose also that the  prime $p$ is generic for $G$,
%that $V$ is  non-singular $\bmod\, p$, 
and that $G_0^p$ is a $C$-group. Then $G^p$ is a $C$-group.
\end{THM}

\noindent
\textbf {Proof.} The subgroup $G_{n-1}^p$ of $G^p$ leaves invariant the subspace $V_{n-1}$
of $V$.  Since $p$ is generic for $G$, we may conclude from 
\cite[Lemma 3.1]{monsch2} that this action is faithful. Thus $G_{n-1}^p$ is a string $C$-group
of Euclidean type, as described in \cite[\S 6B]{monsch1}. By Proposition~\ref{keyprops}(a), 
our task is therefore to show that 
$G_0^p \cap G_{n-1}^p = G_{0, n-1}^p$; so consider any $ g \in G_0^p \cap G_{n-1}^p$.
Now since $G_{n-1} = T \rtimes G_{0, n-1}$ projects onto $G_{n-1}^p$, 
we can multiply $g$ by a suitable element of $G_{0, n-1}^p$, and thereby assume that 
$g \in T^p$. We want to show that $g = e$.

We observed earlier that  $g$ acts as a transvection on $V_{n-1}$, 
with $g(x) -x \in \langle c\rangle =\mbox{rad}(V_{n-1})$ for all $x \in V_{n-1}$. 
On the other hand, 
since $g \in G_0^p \cap G_{n-1}^p$ we have $g(x)-x \in V_{0,n-1}$. Finally, we observe that 
$ \langle c\rangle \cap V_{0,n-1} = \{0\}$ by direct inspection of the various cases 
exhibited in Table 1,  taking $m = n-1$.   
(In most cases this trivial intersection is   implied directly by the fact that 
$G_{0,n-1}$ is of spherical type.)

Thus  $g(x) = x$ for all $x \in V_{n-1}$. Invariably for us the final node $n-1$ 
in the diagram $\Delta(G)$ will be connected to node $n-2$, so that 
$\mbox{disc}(V) \sim -\mbox{disc}(V_{n-2,n-1})$ by a 
dual version of (\ref{discdet}).
The latter discriminant is non-zero for all groups $G_{n-1}$ encountered here, 
again because $p$ is generic for $G$. 
Since $V_{n-1}$ is therefore a singular subspace of the non-singular space $V$, we 
conclude from \cite[Th. 3.17]{art} that $g = e$. This completes the proof in all
important cases.
(It is possible that nodes $n-1$, $n-2$ be non-adjacent; but then it is easy to check 
directly that $G^p \simeq G_{n-1}^p \times C_2$ is a $C$-group.)
\hfill $\square$

\vspace*{5mm}

\section{The  $3$-infinity groups}
\label{higher}

The large number of crystallographic Coxeter groups 
$G = [p_1, \ldots, p_{n-1}]$ of higher ranks 
makes it difficult to 
fully enumerate the regular polytopes obtained by our method. However, 
it is clear that many 
interesting examples occur. As a test of our methods, 
we survey in this section 
groups 
  	$$ G = [ \ldots, 3^k, \infty^l, 3^m, \ldots]\;,$$
of general rank $n$ and having all periods 
$p_j \in \{3, \infty\}$.

When $p_j = \infty$, it is convenient to employ the basic system 
defined by the subdiagram
$\ldots \frac{}{\;\;}\!\!\stackrel{1}{\bullet} 
=\!\!=\!\!= \stackrel{1}{\bullet}\!\!\frac{}{\;\;}\ldots$
on nodes $j-1 , j$. (Thus, $1 = b_{j-1}^2 = b_j^2 = -b_{j-1}\cdot b_j$.)
Typically then, $\Delta(G)$ consists of alternating strings 
of single and doubled branches, as in
$$
\ldots \frac{}{\;\;}\!\!\stackrel{1}{\bullet} 
=\!\!=\!\!= \stackrel{1}{\bullet}
\!-\!\!\!-\!\!\!-\!  \stackrel{1}{\bullet}
\!-\!\!\!-\!\!\!-\!  \stackrel{1}{\bullet}
=\!\!=\!\!= \stackrel{1}{\bullet}
=\!\!=\!\!= \stackrel{1}{\bullet}
=\!\!=\!\!= \stackrel{1}{\bullet}\!\!\frac{}{\;\;}\ldots\; .
$$
For the prime $p=3$, each rotation $r_{j-1} r_j$ in $G^3$ has period $3$, 
and we clearly obtain
$$
G^3 \simeq A_n \simeq S_{n+1}\;,
$$
regardless of the allocation of branches \cite[6.1]{monsch1}.
Likewise, if no $p_j = \infty$, then 
$G^p \simeq A_n$ for any prime $p\geq 3$. Thus, we may henceforth
assume
when it suits us that $p\geq 5$ and that $\Delta(G)$ 
has at least one doubled branch. If in this case $V$ is non-singular, then 
$G^p \simeq O_1(V)$ is of orthogonal type 
(see \cite[Th. 3.1]{monsch1}).

Our approach now must be inductive  on the size of certain classes 
of subdiagrams in $\Delta(G)$; but first we must determine the orthogonal structures
on $V, V_0, V_{n-1}$ and $V_{0,n-1}$.

For $n \geq 1$ we let $d_n := \mbox{disc}(V)$ be the discriminant
of the  underlying  basic system for $G = [\infty^{n-1}]$,
as encoded in the diagram
$$
\Delta(G) \;\; = \;\; \stackrel{1}{\bullet} 
=\!\!=\!\!= \stackrel{1}{\bullet}=\!=\!= \stackrel{1}{\bullet}= 
\ldots =\stackrel{1}{\bullet} =\!\!=\!\!= 
\stackrel{1}{\bullet}=\!\!=\!\!= \stackrel{1}{\bullet}  
$$ 
(on $n$ nodes).
From (\ref{discdet}) we 
have  
$d_n  = 1 d_{n-1} - 1^2 d_{n-2} = d_{n-1} - d_{n-2}$, 
for $n \geq 2$ and taking $d_0 := 1$. Thus
\begin{equation}\label{infdisc}
d_n = \left\{
\begin{array}{rr}
1& \mbox{ if } n \equiv 0,1 \bmod 6,\\
0& \mbox{ if } n \equiv 2,5 \bmod 6,\\
-1&\mbox{ if } n \equiv 3,4 \bmod 6.
\end{array}
\right. 
\end{equation}
Again using (\ref{discdet}), we find that the basic system
underpinning $[3^{n-1}]$ has    
$\mbox{disc}(V) = (n+1)/2^n$.
A routine induction then gives the discriminant
$e_{k,l,m}$ corresponding to the basic system for the
group $G = [3^k , \infty^l , 3^m]$,
with $k+l+m = n-1$ and $k,l,m \geq 0$. Thus
\begin{equation}\label{disceklm}
e_{k,l,m} =\frac{1}{2^{k+m+2}} [\,d_{l+1} (4+2k+2m) - d_{l-1}(2k + 2m + 3km)\,]\;.
\end{equation}

In certain  singular cases we have this 
\begin{LEM}\label{inflem}
Let $G = [3^k, \infty^l]$, with $k + l = n-1$ and $l \geq 1$. 
Suppose that the corresponding space $V$ is singular for the prime $p$. Then
$G^p = \widehat{O}_1(V)$.
\end{LEM}
\noindent\textbf{Proof}. Clearly $G^p \leq \widehat{O}_1(V)$.  Now suppose that
$c = \sum_{j=0}^{n-1} x_j b_j \in \mbox{rad}(V)$. Note that each 
scalar $b_{j-1} \cdot b_j \in \{ -1/2 , -1\}$ and is therefore
invertible in $\mathbb{Z}_p$. Thus,  from $0 = b_0 \cdot c = x_0 +(b_0 \cdot b_1) x_1$
we obtain $x_1 = \alpha_1 x_0$, where $\alpha_1 \in \{1,2\}$. Since $\Delta(G)$
is a tree, we can continue to solve for $x_2, \ldots, x_{n-1}$ as multiples of 
$x_0$ to obtain $x_j = \alpha_j x_0$ for various $\alpha_j$ (with $\alpha_0 := 1$),
where  $\alpha_0, \ldots , \alpha_{n-1}$ are determined only by the 
basic system for $G$. 
In the end, as $V$ is singular, the equation   $0 = b_{n-1} \cdot c$ must be 
redundant, so 
we have $\mbox{rad}(V)  = \langle c \rangle$ (with $x_0 \neq 0$).
Anyway,  we may now take  
$c = 1 b_0 + \alpha_1 b_1 + \ldots +\alpha_{n-1} b_{n-1}$. Then 
$V = \langle c \rangle \perp V_0$, where $V_0$ is non-singular  and
$G_0^p \simeq O_1(V_0)$ (because $ l \geq 1$). We thus 
have 
\begin{equation}\label{gensemd}
\widehat{O}_1(V) = T \rtimes G_0^p\;,
\end{equation}
where $T \simeq \mathbb{Z}_p^{n-1}$ is the abelian group generated by transvections
$t_1, \ldots , t_{n-1}$ satisfying 
$t_j(b_i) = b_i + \delta_{i,j} c$
and $t_j(c) = c$, for $1\leq i, j \leq n-1$. 
Now $r_0 \in \widehat{O}_1(V)$ 
induces an 
isometry on $V/{\mbox{rad}(V)} \simeq V_0$. Let $h \in G_0^p$
be the isometry corresponding to $r_0$ under the natural isomorphism
between $O_1(V/{\mbox{rad}(V)})$ and $O_1(V_0) \simeq G_0^p$.  
A short calculation shows that $t_1 = (h^{-1} r_0 )^{q}\in G^p$,
where $q = 1$ or $(p+1)/2$, according as 
$r_0 r_1$ has period $3$ or $\infty$
(in characteristic $0$).
Finally we  show  inductively that  $t_j \in G^p$ for $1 \leq j \leq n-1$; 
this implies that $G^p = \widehat{O}_1(V)$. 
Fixing $j < n-1$ we may suppose $t_i \in G^p$ for all $1 \leq i \leq j$. 
From \cite[Eq. (9)]{monsch1} we note that
$r_j(b_{j-1}) = b_{j-1} + \alpha b_{j}$,
$r_j(b_{j}) = -b_{j}$ and
$r_j(b_{j+1}) = b_{j+1} + \beta b_{j}$,
where in our case the \textit{Cartan integers}  $\alpha, \beta \in \{ 1,2\}$; otherwise,
for $|k-j|>1$, $r_j(b_k) = b_k$. It is then a routine matter to check that
$$ t_{j-1}^{-\alpha} t_j r_j t_j r_j = t_{j+1}^{\beta}\;.$$
It follows by induction that $t_{j+1}^{\beta} \in G^p$ and hence
$t_{j+1} \in G^p$.  
\hfill$\square$

From \cite[\S 5]{monsch1} and \cite[\S 5]{monsch2}, we already know that
$[3^k, \infty^l]^p$ is a string $C$-group for all $p \geq 3$
and ranks $n \leq 4$ (so that  $k + l \leq 3$).
For example, $[\infty,\infty,\infty]^p$ is the automorphism group of a 
self-dual regular $4$-polytope $\mathcal{P}$ of type $\{p,p,p\}$; when $p=5$ 
we find that $\mathcal{P}$ is isomorphic to the classical star-polytope
$\{5, \frac{5}{2},5\}$ (see \cite[Ch. XIV]{rp} and \cite{grstar}).  
Similarly, $[3,\infty,\infty]^5$ gives back the regular 
star-polytope $\{3,5,\frac{5}{2} \}$.

Let us take stock of  our progress so far.
Keep in mind that whenever $G^p$ is a string $C$-group,
so too is the dually generated group, with $k$ and $m$ interchanged.

\begin{THM} Let $G = [3^k, \infty^l]$, with $k + l = n-1$.\label{3inftpols}
Then for all primes  $p \geq 3$, the group $G^p = [3^k, \infty^l]^p$
is a string $C$-group.
\end{THM}
\noindent\textbf{Proof}. For  $ l = 0$ we have already 
observed that $G^p = [3^{n-1}]^p \simeq S_{n+1}$, the group of the $n$-simplex.
Let us now dispose of the case $l = 1$. 
Since the facet group $[3^{k}]$
is spherical, an induction on $k$, together with Theorem~\ref{intersphere}, shows  that
$G^p = [3^k, \infty]^p$ is a string $C$-group.

Now we may suppose $l \geq 2$. From (\ref{disceklm}), with $m = 0$, we have
\begin{equation}\label{ekl0}
e_{k,l,0} = \frac{1}{2^{k+1}} [ (k+2) d_{l+1} - k d_{l-1}]\;.
\end{equation}
We may also suppose that $p \geq 5$. Then it is easy to check that
if any one of the spaces $V, V_0, V_{n-1}, V_{0,n-1}$ is singular
for a given prime $p$, all the others must be non-singular. We also know, as a basis 
for induction,  that 
$[3^k, \infty^l]^p$ is a string $C$-group  whenever $ k+l \leq 3$. Thus we may suppose that
$G_0^p$ and $G_{n-1}^p$ are string $C$-groups. If  $V_0, V_{n-1},
V_{0,n-1}$ are   non-singular, then  the  corresponding subgroups of $G^p$
are all of type $O_1$, since $l \geq 2$. By Theorem~\ref{subspcrit}(a)--(iii) we have 
$G_0^p \cap G_{n-1}^p = O_1(V_{0,n-1}) = G_{0,n-1}^p$. Thus $G^p$ 
is a string $C$-group by Proposition~\ref{keyprops}(a).
If $V_0$ or $V_{n-1}$ is singular, we similarly apply Theorem~\ref{subspcrit}(c).
Finally, if just $V_{0,n-1}$ is singular,  we
employ Theorem~\ref{subspcrit}(b)--(iii)
and apply Lemma~\ref{inflem} to the subspace $V_{0,n-1}$.
\hfill$\square$
% An example from Egon
%

As an example, consider the group $G=[3,3,3,\infty]$ 
of rank $5$. 
Here we obtain regular $5$-polytopes of type $\{3,3,3,p\}$ with group 
$O_{1}(5,p,0)$. A particularly interesting case occurs when $p=5$. 
Then the polytope can be viewed as a regular tessellation of type 
$\{3,3,3,5\}$ on a hyperbolic $4$-manifold 
whose $78000$ tiles (facets) are 
$4$-simplices and whose $650$ vertex-figures are $600$-cells 
(see \cite[6J]{arp}).  
The dual tessellation has 
$120$-cells as tiles and $4$-simplices as vertex-figures. 
(As an aside, when $p=5$, the group $G=[4,3,3,\infty]$ 
similarly gives a regular tessellation of type $\{4,3,3,5\}$ 
with group $O(5,p,0)$ on a 
hyperbolic $4$-manifold, whose 
facets are $4$-cubes and whose vertex-figures are $600$-cells. 
The dual tessellation again has $120$-cells as tiles and 
$4$-crosspolytopes as vertex-figures.) 

Now let us generalize a little and consider $G = [3^k , \infty^l , 3^m]$, with $k + l + m = n-1$. 
All cases with $k l m = 0$ are covered by Theorem~\ref{3inftpols}, so we assume $k,l,m \geq 1$. 
Then the only new string $C$-group   of  small rank is $[3,\infty,3]$.
Again, we tackle $G^p$ inductively; but since the details are more complicated, 
we shall settle for slightly less comprehensive results.

\begin{THM}\label{3inf3pols}
Let $G = [3^k , \infty^l , 3^m]$, with $k + l + m = n-1$ and 
$k,l,m \geq 1$; and suppose $p$ is an odd prime. Then 

{\rm (a)} $G^p = [3, \infty^l, 3]^p$, with $n= l+3, l\geq 1$ is a string $C$-group,
except possibly when $p = 7$ and $l \geq 4$, with $l \equiv 1 \bmod 3$.

{\rm (b)}  $G^p = [3^k, \infty^l, 3^m]^p$, with $k>1$ or $m >1$, and $l\geq 1$,  
is a string $C$-group 
for all but finitely many primes $p$.
\end{THM}

\noindent\textbf{Proof}. The details in  part (a) 
are very similar to those for Theorem~\ref{3inftpols},
which also serves as a basis for our induction.
The  analysis there fails  only when 
the subspaces  $V$, $V_0$, $V_{n-1}$, $V_{0,n-1}$ are singular in
more than one of dimensions $n, n-1$ and $n-2$; this can only 
occur when $l \equiv 1 \bmod 3$ and $p=7$ (in this case exactly 
$V$ and $V_{0,n-1}$ are singular). In fact, when $p = 7$, $l=1$,
we can use GAP   to verify that $G^7$ is a string $C$-group anyway
\cite[p. 347]{monsch1}.
 
For part (b), we notice in (\ref{disceklm}) that
 $e_{k,l,m} = 0$ (in characteristic $0$) only
when $k = m = 0$ and $l\equiv 1 \bmod 3$. Typically then, $G^p$ falls under the
fully non-singular case described  in Theorem~\ref{subspcrit} (a)--(iii).
\hfill$\square$

\vspace*{5mm}
\noindent\textbf{Remarks}. It may well be in the previous Theorem  that $G^p$  is 
a $C$-group for all primes. Certainly in specific cases, one can explicitly list
the `doubtful' primes; but there seems to be little served by trying to do more here.

\vspace*{5mm}

Now we hunt for  contrary cases in which $G^p$ is \textit{not}
a string $C$-group.  Once again we may assume $p \geq 5$.
Evidently we should examine the groups
$G = [\infty , 3^k , \infty]$, with $k \geq 1$, for which we have
$\mbox{disc}(V) = (k-2)/2^{k+1}$.
By Theorem~\ref{3inftpols}, both $G_0^p$ and $G_{n-1}^p$
are string $C$-groups; and 
$\mbox{disc}(V_0) = \mbox{disc}(V_{n-1}) = -k/2^{k+1}$. Finally,
$G_{0, n-1}^p = [3^k]^p$ is the symmetric group of order $(k+2)!$, and
$\mbox{disc}(V_{0, n-1}) = (k+2)/2^{k+1}$. 

\begin{LEM}\label{3kV}
Suppose $G = [3^k]$, acting as usual on the ($k+1$)--dimensional 
space $V$; and let $p \geq 5$.

{\rm (a)} If $ p \nmid (k+2)$, then $V$ is non-singular and
$G^p \leq O_1(V)$, with equality only when $O(V) = O(2,5,-1)$ or $O(2,7,+1)$.

{\rm (b)} If $ p \mid (k+2)$, then $V$ is  singular and
$G^p$ is always a proper subgroup of $\widehat{O}_1(V)$.
\end{LEM}

\noindent\textbf{Proof}. Recall that  $G^p \leq O_1(V), \widehat{O}_1(V)$
in (a), (b) respectively. We have equality in (a) only when $|O_1(V)| = (k+2)!$.
But the highest power of $p$ dividing $(k+2)!$ is $p^\nu$, where
$$ \nu = \lfloor \frac{k+2}{p}\rfloor  + \lfloor \frac{k+2}{p^2}\rfloor + \ldots\,
<\, \frac{k+2}{p} (1 + \frac{1}{p} + \dots)\, =\, \frac{k+2}{p-1}\;,$$
(see \cite[Prop. 2.3.2]{grove}).
Taking $n = k+1$ in \cite[\S 3.1]{monsch1}, we find that the highest power
$p^\mu$ dividing $|O_1(V)|$ for a non-singular space $V$
has $\mu = \lfloor \frac{k^2}{4}\rfloor$.
(Since $p>3$ we again ignore the groups in (\ref{badgps}).)
Clearly, we usually have $\nu < \mu$, and it is easy by inspection to determine the
two cases with  $G^p = O_1(V)$. 
The situation for (b) is similar. The splitting in (\ref{gensemd}) continues to 
hold in the present context (in which $l = 0$); from this  we get 
$\mu = k + \lfloor \frac{(k-1)^2}{4}\rfloor$ and ultimately
no cases of equality at all.
\hfill$\square$

Now we can show that $G^p$ is very often \textit{not} a string $C$-group:

\begin{THM}\label{notCgps}
 Suppose that $G$ has a string subgroup of the form 
 $[\ldots , \infty, 3^k, \infty, \ldots]$,
 with $k \geq 1$. Let $ p \geq 5$.
 Then $G^p$ is not a string $C$-group, except possibly when $p=5$ or $7$ and
 $k \leq 1$ for all such string subgroups $[\ldots , \infty, 3^k, \infty, \ldots]$.
\end{THM}
\noindent
\textbf{Proof}. Since our intention is to show that $G^p$ fails to be 
a string $C$-group, we may assume by Proposition~\ref{keyprops}(b)
that $G = [\infty, 3^k, \infty]$. 

If $V_{0,n-1}$ is singular, then $p \mid (k+2)$, so that $p \nmid k(k-2)$, 
implying that
$V,V_0, V_{n-1}$ are non-singular. By Theorem~\ref{subspcrit}(b)--(iii) and 
Lemma~\ref{3kV}(b) (applied to $V_{0,n-1}$),
we have $G_0^p \cap G_{n-1}^p = \widehat{O}_1(V_{0,n-1}) \neq G_{0,n-1}^p$. Thus
$G^p$ is not a string $C$-group. 

Suppose that $V_{0,n-1}$ is non-singular, with $p>7$ when $k = 1$. Thus, $ p \nmid (k+2)$,
and $G_{0,n-1}^p \neq O_1(V_{0,n-1})$ by Lemma~\ref{3kV}(a).  If 
also $ p \nmid k$, then we have $G_0^p \cap G_{n-1}^p =  O_1(V_{0,n-1})$
by Theorem~\ref{subspcrit}(a)--(iii), so that $G^p$ is not a string $C$-group.

Finally, suppose $p \mid k$. Thus, $k \geq 5$, and this gives enough wiggle 
room to destroy
polytopality in another way. Let $I = \{0,1\}$, $J = \{5, \ldots, n-1\}$.
Recall, for example, that $G^p_J$ denotes the subgroup generated by the complementary
set of reflections $r_0, \ldots, r_4$; these reflections leave invariant the subspace
$V_J$ spanned by $\{b_j\, |\, j \not\in J\} = \{b_0, \ldots,b_4\}$.
But $\mbox{disc}(V_I) \sim (1-k)/2^k$,  $\mbox{disc}(V_J) \sim -3$
and  $\mbox{disc}(V_{I \cup J}) \sim 1/2$, so that these 
subspaces are all non-singular subspaces
of the non-singular space $V$.  Hence, 
$G^p_I = O_1(V_I) \geq O_1(V_{I \cup J})$
and
$G^p_J = O_1(V_J) \geq O_1(V_{I \cup J})$. 
Thus, if $G^p$ is a string $C$-group, we must have
$$O_1(V_{I \cup J}) \leq G_I^p \cap G_J^p = G^p_{I \cup J} = 
\langle r_2, r_3, r_4 \rangle^p \simeq S_4\;.$$
But $O_1(3,p,0)$ has order $p(p^2-1) > 24$, for $p \geq 5$.
\hfill$\square$

We can summarize the results of this section as follows: if
$G = [p_1, \ldots, p_{n-1}]$ has each period $p_j \in \{3, \infty\}$,
then except for a few small primes, we cannot expect $G^p$ to be a string $C$-group if 
two of the $p_j$'s are $\infty$'s separated by a string of $3$'s. That is, only the groups
$[3^k,\infty^l,3^m]$ can give a string $C$-group, and they do for most
primes $p$.
%%%%%%%%%%%%%%%%%%%%%%%%%%%%%%%%%%%%%%%%%%%%%%%%%%%%%%%%%%%
\vspace*{5mm}

\section{Locally toroidal polytopes of ranks $5$ or $6$}
\label{loctor}
The crystallographic string Coxeter groups  of 
spherical or Euclidean type, along
with the associated modular polytopes, were described 
in \cite[\S 5-6]{monsch1}. When the group is spherical with connected diagram
on $m+1$ nodes, we obtain
(up to isomorphism)
familar  convex regular ($m+1$)-polytopes. After central projection, such polytopes  
can usefully  be viewed as regular spherical tessellations of the circumsphere $\mathbb{S}^{m}$.

Likewise,  each Euclidean group $E$ acts
as the full symmetry group of a certain regular tessellation of Euclidean  
space
 $\mathbb{A}^m$. Indeed, $E$ must be one of the Coxeter groups displayed
in Table~\ref{EucCoxGps}, though perhaps with generators specified in dual 
order. 
A \textit{regular} ($m+1$)-\textit{toroid} $\mathcal{P}$ is the quotient of 
such a tessellation by a non-trivial normal subgroup $L$ of translations in $E$. 
Thus every toroid can be viewed 
as finite, regular tessellation of the $m$-torus. We refer to 
\cite[1D and 6D-E]{arp} for a complete classification; briefly, for each group $E$
the distinct toroids are indexed by a \textit{type vector}
${\bf q} := (q^k, 0^{m-k}) = (q,\ldots,q,0,\ldots,0)$, 
where $q \geq 2$ and $k = 1, 2$ or $m$. 
(For $G = [3,3,4,3]$, the case $k=4$ is subsumed by the case $k=1$.)
Anyway, 
$L$ is generated (as a normal subgroup of $G$) by the translation
$$ \overline{t} := t_1^q \cdots t_k^q\;,$$
where $\{t_1, \ldots , t_m\}$ is a standard set of generators for the 
full group $T$ of translations in $E$.
The modular toroids $\mathcal{P}(E^p)$ described in 
\cite[\S 6B]{monsch1} are special instances; with one exception, we had there
${\bf q} = (p,0,\ldots, 0)$.

In this Section, we consider \textit{locally toroidal}  regular polytopes, 
that is, polytopes
in which the facets and vertex-figures are globally spherical or toroidal, 
as described above
(with at least one kind toroidal). 
The $n$-polytopes of this kind have not yet been fully classified, although 
quite a lot has been discovered since Gr\"{u}nbaum first proposed the
problem in the 1970's (see
\cite{Gru}). What is known rests on a 
broad range of ideas, including frequent use of
 unitary reflection groups, and `twisting' 
and `mixing' operations on  presentations for string $C$-groups. 
We refer to  \cite{McS2,McS3,McS4} for some of the original
investigations, or to
\cite[Chs. 7-12]{arp} for a detailed survey of the project. 

As usual, we begin our own investigation
with a crystallographic linear Coxeter group 
$G = \langle r_0, \ldots , r_{m} \rangle$, but 
immediately  
discard degenerate cases in which the underlying
diagram $\Delta(G)$ is  disconnected. (In such cases
$G^p$ is reducible; 
and $\mathcal{P}(G^p)$ has the 
sort of `flatness' described in
\cite[4E]{arp}.)

In \cite{monsch2} we  discussed all locally toroidal $4$-polytopes 
$\mathcal{P}(G^p)$ which arise from our construction. Turning to higher rank 
$n>4$, we observe that any spherical
facet, or vertex-figure, must  be of type $\{3^{n-2}\}$, $\{4, 3^{n-3}\}$, 
$\{ 3^{n-3},4\}$
or $\{3,4,3\}$ ($n=5$ only). Likewise,  the required Euclidean section must 
have type 
$\{4, 3^{n-4}, 4\}$ or when $n=6$, $\{3,3,4,3\}$ or $\{3,4,3,3\}$.
As described in \cite[Lemma 10A1]{arp}, these constraints severely limit 
the possibilities:
in rank $5$, we have just  $G = [4,3,4,3]$ acting on hyperbolic space 
$\mathbb{H}^4$; and
in rank $6$ we have $G = [4,3,3,4,3], [3,4,3,3,3]$ or $[3,3,4,3,3]$,
all acting on $\mathbb{H}^5$.
Thus we may complete our discussion by examining the modular polytopes which
result from these groups in ranks $5$ or $6$.
\medskip

\subsection{Rank $5$: the group $G = [4,3,4,3]$ }

We may suppose that  $G$ has  diagram
$$\stackrel{2}{\bullet}\!\frac{}{\;\;\;\;\;\;}\!
	\stackrel{1}{\bullet}\!\frac{}{\;\;\;\;\;\;}\!
        \stackrel{1}{\bullet}\!\frac{}{\;\;\;\;\;\;}\!
	\stackrel{2}{\bullet}\!\frac{}{\;\;\;\;\;\;}\!
	\stackrel{2}{\bullet} \;.$$ 
Note that $G_0^p \simeq F_4$ is a $C$-group by \cite[6.3]{monsch1}.
It follows at once from  Theorem ~\ref{eucgen}, and a look at Table 1,
that $G^p$ is a $C$-group for any prime $p \geq 3$. 
(Alternatively, 
since the vertex-figures are spherical, we can appeal to Theorem~\ref{intersphere}.)
Thus $\mathcal{P}(G^p)$  is a locally toroidal regular polytope
of type $\{4,3,4,3\}$. 
Its vertex-figures are copies of the $24$-cell $\{3,4,3\}$.
The facets are toroids  $\{4,3,4\}_{(p,0,0)}$, which one 
could construct by identifying opposite square faces of a $p \times p \times p$ cube. 
(See \cite[6.4]{monsch1}; we note that the facet 
and vertex number mentioned
 there should 
 be $p^{n-1}$ rather than $p^n$.) 
Of course, by flipping the diagram end-for-end, we just as easily obtain the dual polytope
of type $\{3,4,3,4\}$.

In order to identify $G^p$ we consult the list of
irreducible reflection groups in \cite[Table 1]{monsch1}. Since $G^p$ 
has an abelian subgroup of order $p^3$ (generated by translations in the facet),
we immediately rule out the long shots $A_5^p$, $B_5^p$ and $D_5^p$ of
orders $6!$, $2^5 5!$ and $2^4 5!$, respectively. 
Note also that the node label $2$ is a square $\bmod p$ if and only if
$p \equiv \pm 1  \pmod{8}$.
Our hand is now forced:
 $\mathcal{P} = \mathcal{P}(G^p)$
has automorphism group
\begin{equation}\label{loctor5}
\Gamma({\cal P}) =
\left\{ \begin{array}{ll}
O_1(5,p,0)\;,  & \mbox{ if }  p \equiv \pm 1  \pmod{8}\\
O(5,p,0)\;,  & \mbox{ if }  p \equiv \pm 3 \pmod{8}
\end{array}\right.
\end{equation}
For any prime $p \geq 3$, $O_1(5,p,0)$ has order 
$p^4 (p^4 -1) (p^2 - 1) $ and index two in $O(5,p,0)$ (see \cite[pp. 300-301]{monsch1}).

The universal locally toroidal polytopes of rank $5$ are completely
described in \cite[12B]{arp}. 
All but three are infinite.
(The three finite instances have facets with type vector 
$(2,0,0)$, $(2,2,0)$ or $(2,2,2)$; but clearly these do not occur when 
 the modulus $p$ is an odd prime.)
In short,  $\mathcal{P}(G^p)$, being finite, 
is never universal for its type.

Using the fact that  $\mbox{disc}(V)  \sim -2$, we compute that 
the central isometry $-e \in G^p$,  except when  $p \equiv   -1  \pmod{8}$.
When $-e \in G^p$, the polytope $\mathcal{P}(G^p)$ doubly covers the quotient
polytope
$\mathcal{P}(G^p/\{\pm e\})$. The latter polytope still has the same
facets and vertex-figures as $\mathcal{P}(G^p)$. To verify these claims, 
we apply  \cite[2E19]{arp}, 
so  we must show that 
$$ \{ \pm e\} \cap G_{n-1}^p G_0^p = \{ e \}\;.$$
Suppose, on the contrary, that  $-e \in G_{n-1}^p G_0^p$. Since 
$G_{n-1}^p = G_4^p = T \rtimes G_{0,4}^p$,
we can assume $ -e = th$, 
where $h \in G_0^p \simeq F_4$
and   $t$ is some transvection. Thus $h = -t^{-1}$, which has period $2p$. 
This is already 
impossible if $p>3$, since $|F_4| = 3\cdot 2^4 \cdot 4!$. Even when $p=3$ 
it is easy to verify the 
contradiction directly.

\vspace*{5mm}
\subsection{Rank $6$: the groups $[3,4,3,3,3]$, $[3,3,4,3,3]$ and 
$[4,3,3,4,3]$}

In rank $6$ we must consider three closely related groups, beginning with
$$ G  = \langle r_0, r_1, r_2, r_3, r_4, r_5\rangle 
\simeq [3,4,3,3,3]\;.$$
 We may describe a basic system (of roots) for $G$ by the diagram

\begin{equation}\label{loctor6A}
\stackrel{1}{\bullet}\!\frac{}{\;\;\;\;\;\;}\!
	\stackrel{1}{\bullet}\!\frac{}{\;\;\;\;\;\;}\!
	\stackrel{2}{\bullet}\!\frac{}{\;\;\;\;\;\;}\!  
	\stackrel{2}{\bullet}\!\frac{}{\;\;\;\;\;\;}\!
	\stackrel{2}{\bullet}\!\frac{}{\;\;\;\;\;\;}\!
	\stackrel{2}{\bullet}\;. 
\end{equation}
Now the subgroup $H = \langle s_0, \ldots , s_5 \rangle$ generated by the reflections 
\begin{equation}\label{Hgens}
(s_0, s_1, s_2, s_3, s_4, s_5) := (r_1, r_0, r_2 r_1 r_2, r_3, r_4, r_5)\end{equation} 
has index $5$ in $G$ and
is isomorphic to $[3,3,4,3,3]$. The basic system of roots attached to  the $s_j$'s
provides  the diagram 
\begin{equation}\label{loctor6B}
\stackrel{1}{\bullet}\!\frac{}{\;\;\;\;\;\;}\!
	\stackrel{1}{\bullet}\!\frac{}{\;\;\;\;\;\;}\!
	\stackrel{1}{\bullet}\!\frac{}{\;\;\;\;\;\;}\!  
	\stackrel{2}{\bullet}\!\frac{}{\;\;\;\;\;\;}\!
	\stackrel{2}{\bullet}\!\frac{}{\;\;\;\;\;\;}\!
	\stackrel{2}{\bullet} 
\end{equation}
for $H$. Another subgroup
$K = \langle t_0, \ldots , t_5 \rangle$ generated by
\begin{equation}\label{Kgens}
(t_0, t_1, t_2, t_3, t_4, t_5) := (r_2, r_1, r_0, r_3 r_2 r_1 r_2 r_3, r_4, r_5)\end{equation}
has index $10$ in $G$, is isomorphic to $[4,3,3,4,3]$ and has diagram 
\begin{equation}\label{loctor6C}
\stackrel{2}{\bullet}\!\frac{}{\;\;\;\;\;\;}\!
	\stackrel{1}{\bullet}\!\frac{}{\;\;\;\;\;\;}\!
	\stackrel{1}{\bullet}\!\frac{}{\;\;\;\;\;\;}\!  
	\stackrel{1}{\bullet}\!\frac{}{\;\;\;\;\;\;}\!
	\stackrel{2}{\bullet}\!\frac{}{\;\;\;\;\;\;}\!
	\stackrel{2}{\bullet}\;. 
\end{equation}
(See \cite[12A2]{arp}.  Each group acts on $\mathbb{H}^5$ with a 
simplicial fundamental domain of finite volume. In \cite{sqf}, 
these indices were computed by  
dissecting   a simplex for $H$ (or $K$)  into copies of the simplex for $G$.)

Let us now survey the three families of locally toroidal polytopes arising from 
reducing these groups $\bmod \, p$. Since $H_0^p$ and $K_0^p$ are $C$-groups by 
\cite[6B]{monsch1}, we conclude from  Theorem~\ref{eucgen} that 
$H^p$ and $K^p$ are string $C$-groups. Similarly, $G^p$ is a string $C$-group 
because of Theorem~\ref{intersphere}.
 
In each case the underlying space $V$ is non-singular for any prime $p \geq 3$ and has
$$\mbox{disc}(V) = 2 \cdot 0 - (-1)^2 \, (1/4) \sim -1\;. $$
Thus  $\epsilon = +1$ (indicating that the Witt index is $3$). 
As in   rank $5$,
we conclude,  for \textit{each} of the three types, that the 
polytope $\mathcal{P}$ has 
automorphism group
\begin{equation}\label{loctor6gp}
\Gamma({\cal P}) =
\left\{ \begin{array}{ll}
O_1(6,p,+1)\;,  & \mbox{ if }  p \equiv \pm 1  \pmod{8}\\
O(6,p,+1)\;,  & \mbox{ if }  p \equiv \pm 3 \pmod{8}
\end{array}\right.
\end{equation}
Of course, this group is differently generated in the three cases, as indicated above.
We observe that the indices $5$ and $10$ in characteristic $0$
must collapse to $1$ under reduction 
$\bmod \, p$.

For any prime $p \geq 3$, $O_1(6,p,+1)$ has order 
$p^6 (p^4 -1)(p^3 - 1)(p^2 - 1) $ and index 
two in $O(6,p,+1)$ (see \cite[pp. 300-301]{monsch1}). 
We have $-e \in \Gamma(\mathcal{P})$   \textit{except when}  
$ p \equiv -1 \pmod{8}$. 

\noindent
\textbf{Remark}. Some of the results established below  were already announced in 
\cite[12C,D,E]{arp}.

\noindent
\textbf{The Polytopes $\mathcal{P} = \mathcal{P}(G^p)$}.

By \cite[6.2]{monsch1} the vertex-figures of $\mathcal{P}(G^p)$ are $5$-cubes
$\{4,3,3,3\}$. The facets of $\mathcal{P}(G^p)$ are toroids
$\{3,4,3,3\}_{(p,0,0,0)}$, each with $3 p^4$ vertices  
(which corrects
the $3 p^n$ mentioned in \cite[6.5]{monsch1}).

The universal polytope covering $\mathcal{P}(G^p)$ is

$$ \mathcal{U}_{G^p} := \{ \,\{3,4,3,3\}_{(p,0,0,0)} \, , \, \{4,3,3,3\}\, \}\;.$$
Recall that  
$\mathcal{U}_{G^p}$   is conjectured
to be finite only when $p=3$ (see the discussion in \cite[12C1]{arp}, 
restricted to the class of polytopes
under consideration here).
As for $p=3$, it is known that 
$  \mathcal{U}_{G^3}$ 
has automorphism group $\Gamma(\mathcal{U}_{G^3})$ of order 
$3 |G^3| = 3 |O(6,3,+1)| = 72 783 360$. Using GAP we find that $\Gamma(\mathcal{U}_{G^3})$
is a split extension of the additive
cyclic group $(\mathbb{Z}_3, +)$ by $G^3$, that is,
\begin{equation}\label{structUG3} 
\Gamma(\mathcal{U}_{G^3}) =\mathbb{Z}_3 \rtimes O(6,3,+1)\;.
\end{equation}
To check this
directly
we exploit the spinor norm, which for $p=3$ we may view as a homomorphism
$\theta: G^3 \rightarrow \{\pm 1\} = \mathbb{Z}_3^*$. Using this we define
an action of $G^3$ on $\mathbb{Z}_3$ by $g z := \theta(g)\det(g) z$, for 
$g \in G^3, \; z \in \mathbb{Z}_3$.
We obtain the semidirect product  
$ \Lambda :=\mathbb{Z}_3 \rtimes G^3$, with identity $(0,e)$ and 
$$ (y,g)(z,h)  =(y+ g z, gh)\;, $$
for all $y,z \in \mathbb{Z}_3$, $ g,h \in G^3$. Note that 
$r_i z = \eta_i z $, where

$$\eta_i =
\left\{ \begin{array}{cc}
-1, & i \leq 1\\
+1, & i >1\;.
\end{array}
\right. $$

Now let $\rho_0:=(1,r_0)$ and $\rho_i:=(0,r_i)$, for $i\geq 1$;
in brief, $\rho_i = (\delta_{i 0}, r_i)$. It is then a routine matter to check 
that the $\rho_i$'s satisfy the standard relations for the Coxeter group 
$[3,4,3,3,3]$. Indeed,

$$
\rho_i \rho_j = (\delta_{i 0}, r_i)(\delta_{j 0}, r_j) =
(\delta_{i 0} + \eta_i \delta_{j 0},r_i r_j)\;,
$$
so that $\rho_i^2 = (0,e)$, for all $i$. Next we get
$$
(\rho_i \rho_j)^2 = (\delta_{i 0}(1+\eta_i\eta_j)+\delta_{j 0}(\eta_i +\eta_j),  
(r_i r_j)^2)\;,
$$
so that $(\rho_i \rho_j)^2 = (0,e)$, whenever $  i< j-1$. Similarly, 
$\rho_{i-1}\rho_i$ has the same period as $r_{i-1} r_i$ for each $i$.
In particular, we note that
$(\rho_0 \rho_1)^3 = (0,e)$, since $3 = 0$.
Last of all we must verify
the required extra  relation for the facet, namely
$$(\rho_4 \sigma \tau \sigma)^3 = (0,e)\;,$$
where $\sigma := \rho_3 \rho_2 \rho_1 \rho_2 \rho_3$, 
$\tau: = \rho_0\rho_1\rho_2\rho_1\rho_0$;
see \cite[6.5]{monsch1}. Here  check first that 
$\sigma =  (0, s) $ and
$\tau =  (0, t)$, where 
$s:= r_3 r_2 r_1 r_2 r_3$, 
$t:=
r_0 r_1 r_2 r_1 r_0$, then observe that
 $r_4 s t s$ has period $3$ in $G^3$.  

Finally, we
observe that the \textit{Petrie element} $h:= r_0 r_1 r_2 r_3 r_4 r_5$  of $G$ has  
characteristic polynomial $x^6 -x^4 -x^3 -x^2 +1$, so that 
$h^{13} = 6 h^5 +9 h^4 + 6 h^3 +3 h^2 -3h -5e$ and hence $h^{13} = e$  in $G^3$.
The corresponding element of $\Lambda$ is 
$\pi := \rho_0 \rho_1 \rho_2 \rho_3 \rho_4 \rho_5 = (1,h)\;.$
Since $\pi^j = ( j , h^j)$, we have $\pi^{13} = (1,e)$; thus
 the $\rho_i$'s generate $\Lambda$, and $\pi$ has period $39$. 
Since  
$\Gamma(\mathcal{U}_{G^3})$ and $\Lambda$ have equal orders, 
we conclude that   
the two groups are isomorphic (and that $\Lambda$ is a string $C$-group with respect to the 
$\rho_i$'s).
Observe that we obtain the modular polytope $\mathcal{P}(G^3)$ from $\mathcal{U}_{G^3}$
by identitfying vertices separated by $13$ steps along   \textit{Petrie polygons}
of $\mathcal{U}_{G^3}$.

\vspace*{5mm}
\noindent
\textbf{The Polytopes $\mathcal{P} = \mathcal{P}(K^p)$}.

The facets of $\mathcal{P}(K^p)$ are toroids $\{4,3,3,4\}_{(p,0,0,0)}$,
while the vertex-figures are toroids $\{3,3,4,3\}_{(p,0,0,0)}$ of another type.
Consulting \cite[pp. 466-467]{arp}, we note that the universal polytope
$$ \mathcal{U}_{K^p} := \{ \,\{4,3,3,4\}_{(p,0,0,0)} \, , \, \{3,3,4,3\}_{(p,0,0,0)}\, \}$$
is conjectured to exist for all  
primes $p \geq 3$ 
and to be infinite for $p>3$; it is known to be finite for 
$p=3$ (and $p=2$, which again is outside our discussion). 
Our construction of ${\cal P}(K^p)$ establishes the existence 
part of this conjecture for all primes $p\geq 3$. In fact, restricting ourselves to 
$p=3$ for the moment, $\Gamma(\mathcal{U}_{K^3})$ 
also has the same order as 
$\Lambda = \langle \rho_0, \ldots , \rho_5 \rangle$
\cite[Table 12E1]{arp}. Taking a cue from (\ref{Kgens}), we let
$$(\tau_0, \tau_1, \tau_2, \tau_3, \tau_4,\tau_5) := 
(\rho_2, \rho_1, \rho_0, \rho_3 \rho_2 \rho_1 \rho_2 \rho_3, \rho_4, \rho_5)\;.$$
It is routine to check that the $\tau_i$'s  satisfy the defining
relations for $\mathcal{U}_{K^3}$. Next consider the new 
Petrie element 
$$\pi_1 :=  \tau_0 \tau_1 \tau_2 \tau_3 \tau_4 \tau_5 = (-1,h_1) \; ,$$
where $h_1 := t_0 t_1 t_2 t_3 t_4 t_5$. Here
$h_1$ has characteristic polynomial 
$x^6 - x^5 -x^4 -x^2 -x +1$, so that
$$ h_1^{13} = 60 h_1^5  + 48 h_1^4  + 24 h_1^3  + 42 h_1^2  + 15 h_1 - 34 e\;,$$
and hence $h_1^{13} \equiv -e \bmod 3$. Thus 
$\pi_1^{13} = ( -1 , -e)$, and both $\pi_1 $ and $h_1$ have period $26$. 
Moreover, we find that the subgroup generated by the 
$\tau_i$'s contains the crucial element
$$ (\tau_1 \pi_1^{13})^2 = (1,-t_1) (1,-t_1) = (-1,e)\; .$$
Thus, $\Gamma(\mathcal{U}_{K^3}) \simeq \Lambda = \langle \tau_0, \ldots, \tau_5\rangle$.

\vspace*{5mm}
\noindent
\textbf{The Polytopes $\mathcal{P} = \mathcal{P}(H^p)$}.

For each prime $p \geq 3$, the polytope $\mathcal{P}(H^p)$ inherits self-duality
from the hyperbolic tessellation $\{3,3,4,3,3\}$.  To verify this claim, we first use 
(\ref{Hgens}) to establish a basic system of roots  $c_i$ for the $s_i$, namely
$$c_0 := b_1, c_1 :=b_0, c_2 :=r_2(b_1) = b_1+b_2, c_3 :=b_3 , c_4:= b_4, c_5 :=b_5\;.$$
(The Gram matrix $[c_i \cdot c_j]$ for these is encoded in diagram (\ref{loctor6B}).)
Over a suitable extension of the field $\mathbb{Z}_p$, we may now define
an isometry $w$ on $V$ by mapping $c_i \mapsto \alpha_i c_{5-i}$, where
$\alpha_i = 1/\sqrt{2}$ for $i\leq 2$,
$\alpha_i = \sqrt{2}$ for $i>2$.
Thus, $w^2 = e$, $w s_j w = s_{5-j}$, and $w$ induces a
 \textit{polarity} in the polytope $\mathcal{P}(H^p)$.
The facets of $\mathcal{P}(H^p)$ are toroids $\{3,3,4,3\}_{(p,0,0,0)}$; its 
vertex-figures are the dual toroids $\{3,4,3,3\}_{(p,0,0,0)}$.

Consulting \cite[12D3]{arp}, we note that the self-dual universal polytope
$$ \mathcal{U}_{H^p} := \{ \,\{3,3,4,3\}_{(p,0,0,0)} \, , \, \{3,4,3,3\}_{(p,0,0,0)}\, \}$$
is conjectured to exist  
for all primes $p \geq 3$ 
and to be infinite for $p>3$, but is actually known to be finite 
only  for $p=3$ (and $p=2$). Our construction of ${\cal P}(H^p)$ again establishes 
the existence part of the conjecture for all primes $p\geq 3$. Unexpectedly,  
considering our previous look at $ \mathcal{U}_{G^3}$ and  $ \mathcal{U}_{K^3}$, 
we have that 
$ \mathcal{U}_{H^3}$ is a $9$-fold cover of $\mathcal{P}(H^3)$ (\cite[Table 12D1]{arp}).
After constructing $ \mathcal{U}_{H^3}$, we will see that the group $\Lambda$ from above
reappears here as the automorphism group for a \textit{non-self-dual} $3$-fold cover 
of $\mathcal{P}(H^3)$; see (\ref{sit}) below. 

To start the construction we use the automorphism induced  on $H^3$ by 
$w$ to extend the earlier 
action of  $H^3$ ( = $G^3$) on $(\mathbb{Z}_3 , +)$ to an action on 
$ \mathbb{Z}_3 \oplus \mathbb{Z}_3$:
$$ g(y_1 , y_2) := (\theta(w g w) \det(g) y_1 , \theta(g) \det(g) y_2)\, ,$$
for all $ y_1, y_2 \in \mathbb{Z}_3 , g \in H^3$. In the semidirect product
$$\Sigma :=  (\mathbb{Z}_3 \oplus \mathbb{Z}_3) \rtimes H^3\;,$$
with multiplication given by 
\[ (y_{1},y_{2},g)\cdot (z_{1},z_{2},h) =
(y_{1} + \theta(w g w) \det(g) z_{1}, y_{2} + \theta(g) \det(g) z_{2}, gh)\; ,\]
we define 
$$ \sigma_i := (\delta_{i 4} , \delta_{i 1}, s_i)\, , 0\leq i \leq 5\; .$$
It is a straightforward calculation to check that these  $\sigma_i$ satisfy the defining relations for  the automorphism group $\Gamma(\mathcal{U}_{H^3})$. The work is halved by first
noting that the map

\begin{equation}\label{duH3}
\begin{array}{clcc}
\delta : & \Sigma & \rightarrow & \Sigma\\
 & (y_1, y_2, g) & \mapsto & (y_2, y_1, w g w)
\end{array}
\end{equation}
defines an involutory automorphism which transposes each pair $\sigma_i , \sigma_{5-i}$.
(Thus $\delta$ must induce the standard polarity on the self-dual universal polytope
$\mathcal{U}_{H^3}$.) It remains only to check that the $\sigma_i$'s generate $\Sigma$, 
since then $\Sigma$ and $\Gamma(\mathcal{U}_{H^3})$, having equal orders, must be isomorphic.

So consider 
$$ \pi_2 := \sigma_0 \sigma_1 \sigma_2 \sigma_3 \sigma_4 \sigma_5 = (-1, -1, h_2)\, ,$$
where $h_2 := s_0 s_1 s_2 s_3 s_4 s_5$ has period $26$
and satisfies $h_2^{13} = -e$ in $H^3$. Then from
$ \gamma_0 := \sigma_0 \pi_2^{13} = (-1, 1, -s_0)$ and dually
$ \gamma_5 := \sigma_5 \pi_2^{-13} = (1, -1, -s_5)$ we obtain
\begin{equation}\label{eqU3}
\gamma_0^2 = (0, -1, e)\;\; , \;\; \gamma_5^2 = (-1, 0,e)\;,
\end{equation}
so that $\Sigma = \langle \sigma_0, \ldots , \sigma_5 \rangle \simeq \Gamma(\mathcal{U}_{H^3})$.
Furthermore, it is clear that the projection
\begin{equation} 
\begin{array}{clcc}
\varphi : & \Sigma & \rightarrow & \Lambda\\
 & (y_1, y_2, g) & \mapsto & (y_2, g)
\end{array}
\end{equation}
yields yet another set of generators $\overline{\sigma}_i = \varphi(\sigma_i)$
for the group $\Lambda$. Thus 
$\Lambda = \langle \overline{\sigma}_0 , \ldots , \overline{\sigma}_5 \rangle$
is the automorphism group for an intermediate polytope $\mathcal{P}(\Lambda)$,
still of type 
$$\{ \,\{3,3,4,3\}_{(3,0,0,0)} \, , \, \{3,4,3,3\}_{(3,0,0,0)}\, \}\, ,$$
 but now a $3$-fold cover of $\mathcal{P}(H^3)$. From (\ref{eqU3}) we get that
$\varphi(\gamma_0^2)$ and $\varphi(\gamma_5^2)$ have different periods in $\Lambda$,
so that $\mathcal{P}(\Lambda)$ is not self-dual. Evidently the other projection
$\varphi^* : (y_1, y_2,  g) \mapsto (y_1, g)$ yields the automorphism group of the dual polytope
$\mathcal{P}(\Lambda)^*$. 
The situation is summarized here:  

\begin{equation}\label{sit} 
\xymatrix{ & {\mathcal{P}(\Lambda)}\ar[dr]_{3 : 1}&\\
          {\mathcal{U}_{H^3}}\ar[ur]_{3 : 1}^{\varphi} \ar[dr]_{3 : 1}^{\varphi^*} & 
          & {\mathcal{P}(H^3)}\\
           & {\mathcal{P}(\Lambda)}^* \ar[ur]_{3 : 1}& } 
\end{equation}

As interesting as the results in this Section are, it seems that
in order to make further progress with the conjectures in
\cite[\S 12 C,D,E]{arp} concerning locally toroidal polytopes, we
must relax our restriction to a prime modulus $p$
in favour of a more general (composite) modulus $s$. This necessitates a   
somewhat different plan of attack, which we shall pursue in
\cite{monschmod}.

\end{document}